\documentclass[11pt]{amsart}
\usepackage{amsthm,amsbsy,amsfonts,amssymb,amsmath,amscd}

\newcommand{\Q}{\mathbb Q}

\newcommand{\C}{\mathbb C}

\newcommand{\A}{\mathbb A}

\makeatletter
\def\l@section{\@tocline{1}{4pt}{1pc}{}{}}
\def\l@subsection{\@tocline{2}{0pt}{2pc}{5pc}{}}
\makeatother

\begin{document}
\copyrightinfo{~2015}{Dinakar Ramakrishnan}
\title{Recovering Cusp forms on GL$(2)$ from \\ Symmetric Cubes}

\author{Dinakar Ramakrishnan}
\address{253-37 Caltech, Pasadena, CA 91125, USA}				
\email{dinakar@caltech.edu}									

\dedicatory{To Ram Murty, with admiration! \\ Janmadina Shubhecchah !}

\maketitle


\section*{Introduction}						
		
Let $F$ be a number field with adele ring $\A_F$, and let $\pi, \pi'$ be cuspidal automorphic representations of GL$_2(\A_F)$, say with the same central character. If the symmetric squares of $\pi$ and $\pi'$ are isomorphic, we know that $\pi'$ will need to be an abelian, in fact quadratic, twist of $\pi$, which amounts to a multiplicity one statement for SL$(2)$ (\cite{Ra2}). It is of interest to ask if the situation is the same for the symmetric cube transfer (from GL$(2)$ to GL$(4)$) constructed by Kim and Shahidi (cf. \cite{KSh1}). In an earlier paper \cite{Ra8}, dedicated to Freydoon Shahidi, we showed that the answer is in the negative: If $\pi$ is of icosahedral type in a suitable sense (which is meaningful even for $\pi$ without an associated Galois representation), there is a cusp form $\pi^*$ on GL$(2)/F$, which we call the ``conjugate'' of $\pi$, having the same symmetric cube, but which is not an abelian twist of $\pi$. (We also showed there that such a $\pi$ is algebraic when the central character $\omega$ is algebraic, and is moreover rational over $\Q[\sqrt{5}]$ when $\omega=1$; $\pi^*$ is in that case the Galois conjugate of $\pi$ under the non-trivial automorphism of the coefficient field.) In this Note we consider the converse direction and show that for $\pi$ not of solvable polyhedral type, if sym$^3(\pi)\simeq {\rm sym}^3(\pi')$ with $\pi'$ not an abelian twist of $\pi$, then a certain degree $36$ $L$-function has a pole at $s=1$. {\it If} one knew the automorphy of the symmetric fifth power of $\pi$, then this pole would imply that $\pi$ is icosahedral with $\pi'$ twist equivalent to the conjugate $\pi^*$. The situation is simpler if one could associate a Galois representation to $\pi$.

Given a cusp form $\pi$ on GL$(2)/F$, one can define, for every $m \geq 1$, an admissible representation sym$^m(\pi)$ of GL$_{m+1}(\A_F)$, and the principle of functoriality predicts that it is automorphic, which is known (without any hypothesis on $\pi)$ for $m \leq 4$ (\cite{GeJ} for $m=2$, \cite{KSh1} for $m=3$, \cite{K} for $m=4$).
We will say that $\pi$ is {\it solvable polyhedral} if
sym$^m(\pi)$ is Eisensteinian for some $m \leq 4$.

\medskip

\noindent{\bf Theorem A} \, \it Let $F$ be a number field, and $\pi, \pi'$ be cuspidal automorphic representations of GL$(2,\A_F)$, such that $\pi$ is not solvable polyhedral, with central characters $\omega$, $\omega'$ respectively. Suppose we have sym$^3(\pi) \, \simeq {\rm sym}^3(\pi')$. Then $\pi'$ is also not solvable polyhedral.
Also, up to replacing $\pi'$ by a cubic twist, $\omega'=\omega$, sym$^4(\pi') \simeq {\rm sym}^4(\pi)$, and exactly one of the following happens:
\begin{enumerate}
\item[(a)] $\pi' \, \simeq \pi$;
\item[(b)] The functorial product $\Pi:= \pi \boxtimes {\rm sym}^2(\pi')$ is cuspidal on GL$(6)/F$, and the formal Rankin-Selberg $L$-function
$L^S(s, {\rm sym}^5(\pi) \times (\Pi\otimes\omega)^\vee)$ is meromorphic and has a pole at $s=1$, where $S$ is any finite set $S$ of places containing the archimedean and ramified primes.
\end{enumerate}
Furthermore, in case (b), we also have
$$
-{\rm ord}_{s=1}L^S(s, {\rm sym}^5(\pi) \times (\Pi'\otimes\omega)^\vee) \, \geq \, 1,
$$
where
$\Pi': = \, {\rm sym}^2(\pi)\boxtimes \pi'$.
\rm

\medskip

The functorial product $\boxtimes$ used above, also called the Rankin-Selberg product (of automorphic forms), from GL$(k)\times {\rm GL}(r)$ to GL$(kr)$ is known to exist for $(k,r)=(2,2)$ (\cite{Ra2}) and for $(k,r)=(2,3)$ by Kim-Shahidi (\cite{KSh1}).

When $\pi$ is of solvable polyhedral type, there is an associated $2$-dimensional $\C$-representation $\rho$ of the global Weil group $W_F$ with solvable image, and the fibres of the symmetric cube may be evaluated directly without recourse to automorphic $L$-functions.

\medskip

\noindent{\bf Corollary B} \, \it Suppose we are in the situation of (b) above, i.e., with $\pi, \pi'$ having the same symmetric cubes and central character, but not twist equivalent. If sym$^5(\pi)$ is in addition automorphic, then it is isomorphic to
$\Pi=\pi \boxtimes {\rm sym}^2(\pi')$, and also to $\Pi':= \pi' \boxtimes {\rm sym}^2(\pi)$. In particular, it is cuspidal.
In this ``icosahedral'' case, sym$^m(\pi)$ is also automorphic for $m=6,7$, and we have isobaric sum decompositions
$$
{\rm sym}^6(\pi) \, \simeq \, (\pi\boxtimes\pi')\otimes\omega^2 \boxplus {\rm sym}^2(\pi')\otimes\omega^2.
$$
and
$$
{\rm sym}^7(\pi) \, \simeq \, {\rm sym}^2(\pi)\boxtimes\pi'\otimes\omega^2 \, \boxplus \, \pi'\otimes\omega^3.
$$
\rm

\medskip

For holomorphic Hilbert modular newforms of weight $\geq 2$ generating $\pi$, there has been a lot of progress recently on the automorphy of the symmetric fifth power of $\pi$ (and more), due to the striking (independent) works of Dieulefait (\cite{Dieu}) and of Clozel and Thorne (\cite{CTh}). But in this situation one can directly describe the fibres of the symmetric cube (and higher powers) by using the openness of the image of $\rho$ \`a l\`a Serre and Ribet. See \cite{Raj} for an elegant general result for $\ell$-adic representations. For $\pi$ defined by a form of weight one, which has an attached Galois representation $\rho$ of Artin type (cf. \cite{W}, Thm.2.4.1), one has a lot of information on the symmetric powers (see \cite{Wang}).

\medskip

When we are in case (b) of the Theorem A, we will say that $\pi$ is of {\it icosahedral type}. In \cite{Ra8}, we used a closely related condition of being $s$-icosahedral (which is equivalent to the one in part (b) of Theorem A above {\it if} sym$^5(\pi)$ is automorphic), and proved that the finite part $\pi_f$ is then $\Q[\sqrt{5}]$-rational when $\omega=1$, but not rational over $\Q$, and that $\pi'$ is its Galois conjugate $\pi^*$ by the non-trivial automorphism of $E$. The fibre has two elements in this case (up to character twists). For general $\omega$, $\pi$ is only rational over the field generated over $E$ by the values of $\omega$.
This paper may be viewed as a completion of \cite{Ra8}, but can be read independently.

\medskip

Here is some general philosophy, which is not needed for this Note, but underlies the motivation. Langlands conjectures that given any cusp form $\pi$ on GL$(n)/F$, there is an associated irreducible, reductive subgroup $G(\pi)$ of GL$_n(\C)$, generalizing the construction of the Zariski closure of the image of a Galois representation $\rho$, if one may be associated to $\pi$. Such a $G(\pi)$ is expected to exist for any $\pi$, even one which is not algebraic, and for any finite-dimensional representation $r$ of GL$_n(\C)$, the way the restriction of $r$ to $G(\pi)$ decomoses should explicate the analytic behavior of $L(s, \pi; r)$. When $n=2$, there are not too many choices for $G(\pi)$, and if it does not contain SL$_2(\C)$, then its image in PGL$_2(\C)$ will be finite or its connected component must be a torus; the former should correspond to when $\pi$ is dihedral, tetrahedral, octahedral or icosahedral, while the latter to when $\pi$ is dihedral, though of infinite image. The solvable polyhedral case is when one is in the dihedral, tetrahedral or octahedral situation, and here one can give an automorphic definition by the works of Kim and Shahidi (\cite{KSh2}). For our study of the fibres of the symmetric cube transfer, the icosahedral case is the one of interest. For a $\pi$ which is not polyhedral, $G(\pi)$ should contain SL$_2(\C)$ and so all the symmetric power $L$-functions of $\pi$ should be entire. It may be remarked that the article \cite{MuRaj} of Ram Murty with Rajan investigates, under a hypothesis, certain analytic consequences of that case {\it of general type}.

\medskip

This article is dedicated to Ram Murty, a friend whose works I have long read with interest.

\vskip 0.2in

\section{Preliminaries}

\medskip

Let $F$ be a number field. If $\pi_1, \dots, \pi_k$ are isobaric automorphic representations of GL$_{n_1}(\A_F), \dots, {\rm GL}_{n_k}(\A_F)$ respectively, and if $r^1, \dots, r^k$ are polynomial representations of GL$_{n_1}(\C), \dots, {\rm GL}_{n_k}(\C)$, then for any (idele class) character $\mu$ of $F$, we have the associate Langlands $L$-function
$$
L(s, \pi_1, \dots, \pi_k; r^1\otimes \dots \otimes r^k\otimes \mu)\leqno(1.1)
$$
of degree $d=\sum_{j=1}^k {\rm dim}(r^j)$, equipped with an Euler product over the places $v$ of $F$, convergent in a right half plane. Let $S$ denote the finite set of places of $F$ made up of the union of the archimedean places and the finite places where some $\pi_j$ or $\mu$ is ramified. Then for every finite place $v$ outside $S$ of norm $q_v$ and uniformizer $\varpi_v$, there are conjugacy classes $A_v(\pi_{1,v}, \dots, A_v(\pi_{k,v}$ in GL$_{n_1}(\C), \dots, {\rm GL}_{n_k}(\C)$ respectively, such that the $v$-factor of (1.1) is given by
$$
L_v(s, \pi_1, \dots, \pi_k; r^1\otimes \dots \otimes r^k\otimes \mu) \, = \,
{\rm det}\left(I_d - q_v^{-s}\mu(\varpi_v){r^1(A_v(\pi_1))\otimes \dots \otimes r^k(A_v(\pi_k))}\right)^{-1}.\leqno(1.2)
$$
Even at
a ramified (resp. archimedean) place $v$, one can use the local
Langlands correspondence, established for GL$(n)$ by Harris-Taylor and Henniart, to define the corresponding local factor, but we will not need to use it.

We will denote by $r_m$ the standard ($m$-dimensional) representation of GL$_m(\C)$, and by sym$^j$, resp. $\Lambda^j$, the symmetric, resp. alternating, $j$-th power of $r_m$. Fr any Euler product, $L(s) = \prod_v L_v(s)$, and a set $T$ of places of $F$, we will denote by
$L^T(s)$ the incomplete Euler product $\prod_{v \not\in T} L_v(s)$.

\medskip

For $k=2$ and $r^j=r_{n_j}$, the $L$-function (2.1) is called the Rankin-Selberg $L$-function of the pair $(\pi_1, \pi_2)$ (\cite{JPSS}, \cite{Sh3}), which admit meromorphic continuation and a standard functional equation, such that when both $\pi_j$ are cuspidal, there is a pole at $s=1$ iff $n_1=n_2$ and moreover, $\pi_2$ is the contragredient $\pi_1^\vee$ (of $\pi_1$). It is expected (by Langlands' functoriality principle) that there exists an isobaric automorphic form $\pi_1 \boxtimes \pi_2$, called the Rankin-Selberg product or the automorphic tensor product, on GL$(n_1n_2)/F$, whose standard $L$-function agrees with $L(s,\pi_1 \times \pi_2)$, which is a shorthand for $L(s, \pi_1, \pi_2; r_{n_1}\otimes r_{n_2})$. This is known to be true for $n_1=n_2 =2$ (\cite{Ra2}) and for $(n_1, n_2) =(2,3)$ by Kim-Shahidi (\cite{KSh1}).

\medskip

For any cuspidal (and hence isobaric) automorphic form $\pi$ on GL$(2)/F$, a fundamental result we will use is the existence, for $j \leq 4$, of the symmetric $j$-th power transfer to GL$(j+1)/F$, which is classical (due to Gelbart-Jacquet) for $j=2$ (cf. \cite{GeJ}) and for $j=3, 4$ due to Kim-Shahidi (\cite{KSh1}) and Kim (\cite{K}). For $j=2$, one knows (by \cite{GeJ}) that (for cuspidal $\pi$), sym$^2(\pi)$ is cuspidal iff $\pi$ is not {\it dihedral}, i.e., monomial attached to a character of a quadratic extension $K/F$. A beautiful result of Kim and Shahidi (\cite{KSh2}) asserts that for $j=3$, resp. $j=4$, sym$^j(\pi)$ is cuspidal iff $\pi$ is not {\it tetrahedral}, resp. {\it octahedral}, which means sym$^2(\pi)$ is not monomial attached to a character of a cyclic cubic, resp. non-normal cubic, extension $E/F$. We will say that $\pi$ is {\it solvable polyhedral} iff it is dihedral, tetrahedral or octahedral.

\medskip

We will call a cusp form $\pi$ {\it s-icosahedral} iff there is another cusp form $\pi^*$ such that
$$
L^S(s, \pi; {\rm sym}^5) \, = \, L^S(s, {\rm sym}^2(\pi^*) \boxtimes \pi\otimes\omega).
$$
We will at times write, by abuse of notation $L(s, {\rm sym}^m(\pi))$ for any $m$, though sym$^m(\pi)$ is only known to be an admissible representation of GL$_{m+1}(\A_F)$ for general $m$.

\medskip

We will also have occasion to make use of Kim's exterior square functoriality $\eta \mapsto \Lambda^2(\eta)$ (cf. \cite{K}) from GL$(4)/F$ to GL$(6)/F$, such that the standard $L$-function of $\Lambda^2(\eta)$ agrees with $L(s, \eta; \Lambda^2)$.

\medskip

\section{Two Lemmas}

\medskip

Let $\pi, \pi'$ be cuspidal automorphic representations of GL$_2(\A_F)$ which are not abelian twists of each other, and {\it not solvable polyhedral}, with respective central characters $\omega, \omega'$, such that sym$^3(\pi)$ and sym$^3(\pi')$ are isomorphic. Let $S$ denote a finite set of places of $F$ containing the archimedean and ramified places (for $\pi, \pi'$).

\medskip

\noindent{\bf Lemma 2.1} \, \it Suppose $\pi, \pi'$ are as above, and are not twist equivalent. Then
\begin{enumerate}
\item[(a)]The automorphic representations
$$
\pi \boxtimes \pi', \, \Pi:= \pi\boxtimes {\rm sym}^2(\pi'), \, \Pi':= \pi'\boxtimes {\rm sym}^2(\pi)
$$
are all cuspidal. Moreover, $\pi\boxtimes\pi'$ does not admit any non-trivial self-twist.
\item[(b)]When $\omega=\omega'$, $\Pi \, \simeq \, \Pi'$.
\end{enumerate}
\rm

\medskip

{\it Proof}. \, (a) \, Just the fact that $\pi, \pi'$ are not dihedral implies that $\pi\boxtimes\pi'$ is cuspidal unless $\pi'$ is an abelian twist of $\pi$ (cf. \cite{Ra2}), which we have assumed to be not the case. So we have cuspidality in this case. Suppose $\pi \boxtimes \pi'$ is isomorphic to $\pi\boxtimes\pi'\otimes\chi$ for a character $chi$. Since we know the fibres of $(\pi, \pi') \mapsto \pi \boxtimes \pi'$ (from \cite{Ra2}), we see that there must be characters $\chi_1, \chi_2$ such that $\chi=\chi_1, \chi_2$, $\pi\simeq\pi\otimes\chi_1$ and $\pi'\simeq\pi'\otimes\chi_2$. Since $\pi, \pi'$ do not admit any self-twist, for otherwise they will be dihedral, we are forced to have $\chi_1=\chi_2=1$, thus $\chi=1$, and the assertion of part (a) is proved for $\pi\boxtimes\pi'$.

Next we will show the assertions for $\Pi$ and note that by symmetry they also hold for $\Pi'$.
Consider the Rankin-Selberg $L$-function
$$
L^S(s, \Pi\times\Pi^\vee),
$$
which can be rewritten as
$$
L^S(s,(\pi \boxtimes \pi^\vee)\times ({\rm sym}^2(\pi')\boxtimes{\rm sym}^2(\pi')^\vee)),
$$
with isobaric decompositions
$$
\pi \boxtimes \pi = {\rm sym}^2(\pi)\otimes\omega^{-1} \boxplus \underline{1}
$$
and
$$
{\rm sym}^2(\pi')\boxtimes{\rm sym}^2(\pi')^\vee \, \simeq \, {\rm sym}^4(\pi')\otimes{\omega'}^{-2} \boxplus ({\rm sym}^2(\pi')\otimes{\omega'}^{-1}) \boxplus \underline{1},
$$
which may be taken as its definition. (Here $\underline{1}$ denotes the trivial automorphic representation of GL$_1(\A_F)$.)
Thus we have the factorization
$$
L^S(s,  \Pi\times\Pi^\vee) \, = \, L^S_1(s)L^S_2(s),
$$
where
$$
L^S_1(s):= \, L^S(s, {\rm sym}^2(\pi)\times {\rm sym}^4(\pi')\otimes(\omega{\omega'}^2)^{-1})L^S(s,{\rm sym}^2(\pi)\times{\rm sym}^2(\pi')\otimes(\omega\omega')^{-1})L^S(s,{\rm sym}^2(\pi)\otimes{\omega}^1),
$$
and
$$
L^S_2(s):= \, L^S(s, {\rm sym}^4(\pi')\otimes{\omega'}^{-2})L^S(s,{\rm sym}^2(\pi')\otimes{\omega'}^{-1})\zeta_F^S(s).
$$

Note that by Jacquet-Shalika (\cite{JS1}), $\Pi$ is cuspidal iff the (incomplete) Rankin-Selberg $L$-function $L^S(s, \Pi\times\Pi^\vee)$ has a simple pole at $s=1$. So we have to show that $L^S_1(s)L^S_2(s)$ has a pole of order one. Since sym$^m(\eta)$ is cuspidal for $\eta=\pi, \pi'$, we see that $L^S_2(s)$ has a simple pole at $s=1$, and that the only possible pole of $L^S_1(s)$ could come from the factor $L^S(s,{\rm sym}^2(\pi)\times{\rm sym}^2(\pi')\otimes(\omega\omega')^{-1})$, which is usually written as $L^S(s, {\rm Ad}(\pi)\times{\rm Ad}(\pi'))$, where Ad$(\pi)$ is the selfdual adjoint sym$^2(\pi)\otimes\omega^{-1}$. This factor can have a pole iff Ad$(\pi)$ and Ad$(\pi')$ are isomorphic, which by \cite{Ra2} can happen iff $\pi$ is an abelian twist of $\pi'$, which we have assumed to be not the case.
We are now done with proving part (a).

(b) \, Suppose $\omega=\omega'$. To prove that $\pi \simeq \i'$, we only have to show the existence of a pole at $s=1$ of
$$
L^S(s, \Pi \times {\Pi'}^\vee) \, = \, L^S(s, (\pi\boxtimes{\rm sym}^2(\pi)^\vee)\times \left(\pi'\boxtimes{\rm sym}2(\pi')^\vee\right)^\vee).
$$
We have
$$
\pi\boxtimes{\rm sym}^2(\pi)^\vee \, \simeq \, {\rm sym}^3(\pi)\omega^{-2} \boxplus \pi\otimes\omega^{-1},
$$
and similarly for the corresponding expression involving $\pi'$. Thus $L^S(s, \Pi \times {\Pi'}^\vee)$ factors as
$$
L^S(s,{\rm sym}^3(\pi)\times{\rm sym}^3(\pi')^\vee)L^S(s,{\rm sym}^3(\pi)\times{\pi'}^\vee\otimes\omega^{-1})L^S(s,\pi\times{\rm sym}^3(\pi')^\vee\otimes{\omega'}^{-1})
$$
times the entire function (since $\pi, \pi'$ are inequivalent):
$$
L^S(s,\pi\times{\pi'}^\vee).
$$
Note that $L^S(s,{\rm sym}^3(\pi)\times{\rm sym}^3(\pi')^\vee)$ has a simple pole at $s=1$ since the symmetric cubes of $\pi, \pi'$ are cuspidal and equivalent. The remaining two $L$-functions dividing
$L^S(s, \Pi \times {\Pi'}^\vee)$ are entire and non-zero at $s=1$.
\qed

\medskip

\noindent{\bf Lemma 2.2} \, \it Up to replacing $\pi$ by a cubic twist, we have
\begin{enumerate}
\item[(a)]$\omega \, = \, \omega'$.
\item[(b)]${\rm sym}^4(\pi) \, \simeq \, {\rm sym}^4(\pi')$.
\end{enumerate}
\rm

\medskip

{\it Proof}. \, At any $v \not\in S$ with uniformizer $\varpi_v$, let the corresponding conjugacy classes $A_v(\pi), A_v(\pi')$ (of $\pi, \pi'$) be represented by diag$(\alpha_v,\beta_v)$, ${\rm diag}(\alpha_v',\beta_v')$ respectively, so that
$$
\omega_v(\varpi_v)=\alpha_v\beta_v, \, \omega'_v(\varpi_v)=\alpha'_v\beta'_v.
$$
A direct calculation shows
$$
\Lambda^2({\rm sym}^3(\pi_v) \, \simeq \, {\rm sym}^4(\pi_v)\otimes\omega_v \, \oplus \, \omega_v^3.\leqno(2.4)
$$
This is the $v$-factor of the exterior square of the cusp form sym$^3(\pi)$ by Kim \cite{K}. By the strong multiplicity one theorem for isobaric automorphic representations (\cite{JS1}, we obtain an isomorphism
$$
\Lambda^2({\rm sym}^3(\pi)) \, \simeq \, {\rm sym}^4(\pi)\otimes\omega \, \boxplus \, \omega^3.\leqno(2.5)(i)
$$
Similarly,
$$
\Lambda^2({\rm sym}^3(\pi)) \, \simeq \, {\rm sym}^4(\pi')\otimes\omega' \, \boxplus \, {\omega'}^3.\leqno(2.5)(ii)
$$
Now since $\pi, \pi'$ are not solvable polyhedral, sym$^4(\pi)$ and sym$^4(\pi')$ are cuspidal. f $\omega^3$ is distinct from ${\omega'}^3$, then we must have
$$
-{\rm ord}_{s=1} L^S(s, {\rm sym}^4(\pi')\otimes\omega'\omega^{-3}) \, \geq 1,
$$
which is impossible by the cuspidality of sym$^4(\pi')$. Hence $\omega^3={\omega'}^3$, yielding art (a).

Comparing (2.5)(i) and (2.5)(ii), we also get part (b).

\qed

\medskip

\section{Proof of Theorem A}

Let $\pi, \pi'$ be as in Theorem A. Since $\pi$ is not solvable polyhedral, sym$^3(\pi)$ is cusidal and does not admit a quadratic self-twist. As $\pi'$ and $\pi$ have the same symmteir cubes, the same statements result for sym$^3(\pi')$, implying by \cite{KSh2} that $\pi'$ is also not solvable polyhedral.

There is nothing to prove if $\pi$ and $\pi'$ are abelian twists of each other, so we may assume that they are not. By the Key Lemma 2.2, $\omega^3={\omega'}^3$. Thus, if we write $\nu$ for $\omega/\omega'$ and consider $\pi''=\pi'\otimes\nu^{-1}$. Then sym$^3(\pi'') \simeq {\rm sym}^3(\pi')$ and moreover, the central character $\omega''$ of $\pi''$ is $\omega'\nu^{-2}$, which is $\omega'\nu = \omega$. Thus up to replacing $\pi'$ by a cubic twist, we may assume that $\pi$ and $\pi'$ have the same central character.

\medskip

Consider the functorial product $\pi \boxtimes \pi'$ which, by Lemma 2.1, is cuspidal and not isomorphic to any non-trivial abelian twist of itself. We have, for any character $\mu$ of $F$,
$$
L^S(s, {\rm sym}^4(\pi)\times (\pi\boxtimes\pi')\otimes \mu) \, = \, L^S(s, {\rm sym}^5(\pi)\times \pi'\otimes\mu)L^S(s, {\rm sym}^3(\pi)\times \pi'\otimes\mu\omega),\leqno(3.1)
$$
which also equals (by replacing sym$^4(\pi)$ by sym$^4(\pi')$ in the left hand side $L$-function of (3.1))
$$
L^S(s, {\rm sym}^5(\pi')\times \pi\otimes\mu)L^S(s, {\rm sym}^3(\pi')\times \pi\otimes\mu\omega).\leqno(3.2)(a)
$$
We may rewrite (3.2)(a) by replacing sym$^3(\pi')$ in (3.2) by sym$^3(\pi)$ and decomposing
${\rm sym}^3(\pi)\times \pi$, to obtain
$$
L^S(s, {\rm sym}^5(\pi')\times \pi\otimes\mu)L^S(s, {\rm sym}^4(\pi)\otimes\mu\omega)L^S(s, {\rm sym}^2(\pi)\otimes\mu\omega^2).\leqno(3.2)(b)
$$
Appropriately twisting (3.1) and (3.2(b) by sym$^2(\pi)^\vee\otimes(\mu\omega^2)^{-1}$ ($={\rm sym}^2(\pi)\otimes(\mu\omega^4)^{-1}$), we are led to identify
$$
L^S(s, {\rm sym}^5(\pi)\times {\rm sym}^2(\pi)\times \pi'\otimes\omega^{-4})
L^S(s, {\rm sym}^3(\pi)\times {\rm sym}^2(\pi)\times \pi'\otimes\omega^{-3})\leqno(3.3)(a)
$$
with
$$
L^S(s, {\rm sym}^5(\pi')\times {\rm sym}^2(\pi)\times \pi\otimes\omega^{-4})L^S(s, {\rm sym}^4(\pi)\times {\rm sym}^2(\pi)\otimes\omega^{-4})L^S(s, {\rm sym}^2(\pi)\times{\rm sym}^2(\pi)^\vee).\leqno(3.3)(b)
$$
We have, for any character $\nu$ of $F$,
$$
L^S(s, {\rm sym}^5(\pi')\times {\rm sym}^2(\pi)\times \pi\otimes\nu) \, = \, L^S(s, {\rm sym}^5(\pi')\times{\rm sym}^3(\pi')\otimes\nu)L^S(s, {\rm sym}^5(\pi')\times \pi'\otimes\nu\omega),\leqno(3.4)(a)
$$
which equals
$$
L^S(s, {\rm sym}^8(\pi')\otimes\nu)L^S(s, {\rm sym}^6(\pi')\otimes\nu\omega)^2L^S(s, {\rm sym}^4(\pi')\otimes\nu\omega^2).\leqno(3.4)(b)
$$
The proof of Lemma 2.2 shows that the $L$-function product of (3.4)(b) is invertible at $s=1$. Using this in conjunction with (3.4)(a), we see that the expression in (3.3)(b) has a simple pole at $s=1$, which results in a simple pole (at $s=1$) of (3.3)(a). On the other hand,
$$
L^S(s, {\rm sym}^3(\pi)\times {\rm sym}^2(\pi)\times \pi'\otimes\omega^{-4}) \, = \, L^S(s, {\rm sym}^3(\pi')\times {\rm sym}^2(\pi)\times \pi'\otimes\omega^{-4}),\leqno(3.5)(a)
$$
which factors as
$$
L^S(s, {\rm sym}^4(\pi')\times {\rm sym}^2(\pi)\otimes\omega^{-4})L^S(s, {\rm sym}^2(\pi')\times {\rm sym}^2(\pi)\otimes\omega^{-3}).\leqno(3.5)(b)
$$
The $L$-function on the left of (3.5)(b) has no pole since (by virtue of $\pi'$ not being solvable polyhedral) sym$^4(\pi')$ is cuspidal. And the one on the right has a pole iff Ad$(\pi): = {\rm sym}^2(\pi)\otimes\omega^{-1}$ is isomorphic to Ad$(\pi')$, which would imply, by \cite{Ra2}, that $\pi$ and $\pi'$ are twist equivalent, which is not the case. Thus the only possibility, by looking at (3.3)(a), is to have
$$
-{\rm ord}_{s=1} L^S(s, {\rm sym}^5(\pi)\times \Pi\otimes\omega^{-4}) \, = \, 1,(3.6)
$$
where, as before,
$$
\Pi \, = \, {\rm sym}^2(\pi)\boxplus \pi'.
$$
We get the identification of sym$^5(\pi)$ with $\Pi\otimes\omega$, once we note that
$$
(\Pi\otimes\omega)^\vee \, \simeq \Pi\otimes\omega^{-4}.
$$
Since we know by Lemma 2.1 that $\Pi$ and $\Pi'$ are isomorphic, we get the dichotemy of Theorem A.

\qed

\medskip

\section{Proof of Corollary B}

\medskip

Suppose in addition to our working hypotheses, we also know that sym$^5(\pi)$ is automorphic. As before, by replacing $\pi'$ by a cubic twist if necessary, we may assume that it has the same central character $\omega$ as $\pi$.
Then in case (b) of Theorem A, the fact that
$$
-{\rm ord}_{s=1} L^S(s, {\rm sym}^5(\pi) \times (\Pi\otimes\omega)^\vee) \, \geq 1
$$
implies that $\Pi\omega$ occurs in the isobaric decomposition of sym$^5(\pi)$. But both $\Pi$ and sym$^5(\pi)$ are representations of GL$_6(\A_F)$, which forces the isomorphism
$$
{\rm sym}^5(\pi) \, \simeq \, \Pi\otimes\omega = \pi \boxtimes {\rm sym}^2(\pi')\otimes\omega,\leqno(4.1)
$$
which is also isomorphic to $\Pi'\otimes\omega=\pi'\boxtimes{\rm sym}^2(\pi)\omega$.

\medskip

Now we have the identifications
sym$^2(\pi)\boxtimes\pi = {\rm sym}^3(\pi)\boxplus \pi\otimes\omega$, sym$^3(\pi) \simeq {\rm sym}^3(\pi')$,
and
$$
\pi'\boxtimes {\rm sym}^3(\pi')= \, {\rm sym}^4(\pi') \boxplus {\rm sym}^2(\pi')\otimes\omega,\leqno(4.2)
$$
which can be taken to be the definition of $\pi'\boxplus{\rm sym}^3(\pi')$. Thus
we are able to realize, in our (icosahedral) case, the functorial product of sym$^5(\pi)$ and $\pi$ by setting
$$
{\rm sym}^5(\pi)\boxtimes\pi \, \simeq \, {\rm sym}^4(\pi')\otimes\omega \boxplus {\rm sym}^2(\pi')\otimes\omega^2 \boxplus \pi'\boxtimes\pi\otimes\omega^2.\leqno(4.3)
$$

On the other hand, we have the factorization
$$
L(s, {\rm sym}^5(\pi)\times \pi) \, = \, L(s, {\rm sym}^6(\pi))L(s,{\rm sym}^4(\pi)\otimes\omega),\leqno(4.4)
$$
which is in fact correct at every place by the work of Shahidi (\cite{Sh1}). By Lemma 2.2, sym$^4(\pi)$ and sym$^4(\pi')$ are isomorphic, implying (by a comparison of (4.3) and (4.4), that
$$
L(s, {\rm sym}^6(\pi)) \, = \, L(s, {\rm sym}^2(\pi')\otimes\omega^2)L(s, \pi'\boxtimes\pi\otimes\omega^2.\leqno(4.5)
$$
Thus we may realize sym$^6(\pi)$ as the isobaric automorphic representation
$$
{\rm sym}^2(\pi')\boxtimes\pi\otimes\omega^2 \, \boxplus \, \pi'\boxtimes\pi\otimes\omega^2.
$$

The assertion about sym$^7(\pi)$ is established in a similar manner, and the proof is left as an exercise.

\qed

\medskip

\bibliographystyle{math}
\bibliography{fibresSymCube}

\medskip

Dinakar Ramakrishnan

253-37 Caltech

Pasadena, CA 91125, USA.

dinakar@caltech.edu


\end{document}